\documentclass[11pt]{article}
\usepackage{geometry}
\usepackage{graphicx}
\usepackage{setspace}
\usepackage{verbatim}
\usepackage{lmodern}
\usepackage{url}
\usepackage{latexsym}
\usepackage{amsmath,amsthm,amsfonts,amssymb,amscd,amsbsy}
\newcommand{\al}{{\alpha}}
\newcommand{\re}{{\mathbb{R}}}

\newcommand{\R}{\mathbb R}
\newcommand{\N}{\mathbb N}
\newtheorem{theorem}{Theorem}[section]

\newtheorem*{question*}{Question}

\theoremstyle{definition}

\newtheorem{remark}[theorem]{Remark}

\begin{document}

\title{Homoclinic Bifurcations: our collaboration with Jean-Christophe Yoccoz 
  }

\author{Carlos Gustavo Moreira and Jacob Palis \\
IMPA}

\maketitle
\begin{abstract}
We briefly describe our works in collaboration with Jean-Christophe Yoccoz, a great mathematician and friend, with special emphasis on those related to Homoclinic Bifurcations and Fractal Geometry. We also tell some related personal stories. 
\end{abstract}

%\begin{classification}
%Primary 37C29, 28A80
%\end{classification}

%\begin{keywords}
%Fractal geometry. Homoclinic bifurcations. 
%\end{keywords}

%\maketitle

\section{Introduction}
Jacob Palis was visiting the Institut des Hautes \'Etudes Scientifiques-IH\'ES about 1980. One day, Albert Fathi, a French mathematician, who went to IH\'ES to attend a conference, talked to Jacob about an excellent student of Michel Herman, Jean-Christophe Yoccoz, that would, a few hours later, give a talk at nearby Universit\'e de Paris XI at Orsay. Jacob went to hear the talk and was very impressed with it. Immediately after that he invited Jean-Christophe to visit the Institute for Pure and Applied Mathematics (IMPA) at Rio de Janeiro. Jean-Christophe seemed to appreciate the idea and indicated the possibility of spending his French military service in Rio.
He came to Brazil some time after Etienne Ghys, another excellent young French mathematician, who also spent the period of his military service at IMPA.

The first joint works by Jacob and Jean-Christophe were done in the eighties (see \cite{PY1}, \cite{PY2}), and concerned centralizers of diffeomorphisms: They proved that for ``typical" diffeomorphisms in several dynamical contexts, their centralizers are trivial. That is, if $\varphi$ is such a diffeomorphism and $\psi$ is another one such that $\varphi \circ \psi=\psi\circ\varphi$, then $\psi$ is an iterate of $\varphi$: $\psi=\varphi^ k$, for some $k\in\mathbb Z$.  Influenced by that, the master degree thesis of Carlos Gustavo Moreira (Gugu), advised by Jacob in 1990 (when Gugu was 16 years old), was much related to this topic: it was based on a previous work by Nancy Kopell, who proved that centralizers of ``typical" diffeomorphisms of the circle are trivial. The next joint paper by  Jacob and Jean-Christophe (\cite{PY3}), published in 1990, presents a complete set of invariants for differentiable conjugacies of the so-called Morse-Smale diffeomorphisms, an important class of dynamical systems which exist in any compact manifold and are structurally stable, as proved by Jacob in his Ph.D. thesis in low dimensions and by Jacob and Steve Smale in general.

After these works, the collaboration of Jacob (and Gugu) with Jean-Christophe was related to fractal geometry and homoclinic bifurcations. Homoclinic bifurcations are perhaps the most important mechanism that creates complicated dynamical systems from simple ones. This phenomenon takes place when an element of a family of dynamics (diffeomorphisms or flows) presents a hyperbolic periodic point whose stable and unstable manifolds have a non-transverse intersection. When we connect, through a family, a dynamics with no {\it homoclinic points} (namely, intersections of stable and unstable manifolds of a hyperbolic periodic point) to another one with a {\it transverse} homoclinic point (a homoclinic point where the intersection between the stable and unstable manifolds is transverse) by a family of dynamics, we often go through a homoclinic bifurcation. The existence of transverse homoclinic points implies that the dynamics is quite complicated, as remarked already by Poincar\'e: ``Rien n'est plus propre \`a nous donner une id\'ee de la complication du probl\`eme des trois corps et en g\'en\'eral de tous les probl\`emes de Dynamique...", in his classic {\it Les M\'ethodes Nouvelles de la M\'ecanique C\'eleste} (\cite{Po}), written in late 19th century.

Incidentally, Jacob was visiting the Institut Mittag-Leffler at Djursholm, near Stockholm, invited by Lennart Carleson about 1990 and when visiting the Institute's Library, he looked for papers by Poincar\'e, and had the impression of seeing two versions of one of Poincar\'e's article on dynamics. Much later, he became aware of the prize given to Poincar\'e by King Oscar II of Sweden. In fact, there were two versions of Poincar\'e's work presented for the prize - a mistake was detected by the Swedish mathematician Phragm\'en in the first version. When Poincar\'e became aware of that, he rewrote the paper, including the quotation above calling the attention to the great complexity of dynamical problems related to homoclinic intersections.  We refer to the paper \cite{Y} by Jean-Christophe describing such an event.

At this point it is important to mention the notion of {\it hyperbolic systems}, introduced by Smale in the sixties, after a global example provided by Anosov, namely the diffeomorphism $f(x,y)=(2x+y, x+y)\pmod 1$ of the torus ${\mathbb T}^2={\mathbb R}^2/{\mathbb Z}^2$ and a famous example, given by Smale himself, of a {\it horseshoe}, that is a robust example of a dynamical system on the plane with a transverse homoclinic point as above, which implies a rich dynamics - in particular the existence of infinitely many periodic orbits. The figure below depicts a horseshoe.

%\begin{figure}
%\centering
\begin{center}
\includegraphics[width=5cm]{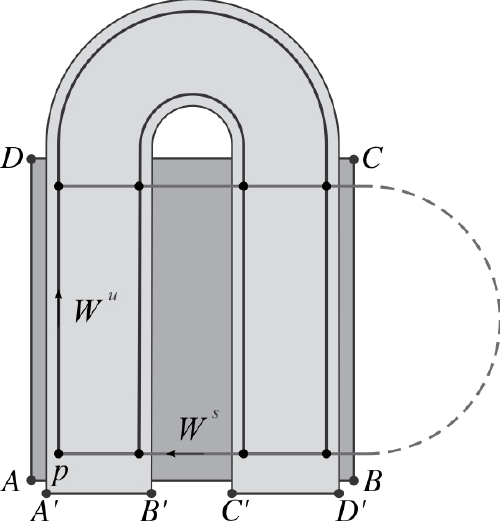}
\end{center}
%\end{figure}

(the dynamics sends the square $(ABCD)$ onto the domain bounded by $(A'B'C'D')$).

Let $\Lambda\subset M$ be a compact subset of a manifold $M$.  We say  that $\Lambda$ is a hyperbolic set for a diffeomorphism $\varphi\colon
M\to M$ if $\varphi(\Lambda) = \Lambda$ and there is a decomposition $T_\Lambda M
= E^s\oplus E^u$ of the tangent bundle of $M$ over $\Lambda$ such that $D\varphi\mid_{E^s}$ is uniformly contracting and
$D\varphi\mid_{E^u}$ is uniformly expanding. We say that $\varphi$ is hyperbolic if the limit set of its dynamics is a hyperbolic set.

Homoclinic bifurcations become important when going beyond the hyperbolic theory. In the late sixties, Sheldon Newhouse combined homoclinic bifurcations with the complexity already available in the hyperbolic theory and some new concepts in Fractal Geometry to obtain dynamical systems far more complicated than the hyperbolic ones. Ultimately this led to his famous result on the coexistence of infinitely many periodic attractors. Later on, Jacob and Floris Takens developped these techniques combining Fractal Geometry and Dynamical Systems in their work on hyperbolicity or lack of it near homoclicic bifurcations, which we shall discuss in the next section. 

Gugu met Jean-Christophe around 1992, introduced by his advisor Jacob, when he was a Ph.D. student at IMPA. His thesis concerned a conjecture by Jacob on regular Cantor sets (see section 3 for a precise definition), which we shall discuss in the next sections. In his thesis ([M]), Gugu made some progress on the question. After his Ph.D., concluded in 1993, he went to Paris for a Postdoctoral program under the guidance of Jean-Christophe at the Universit\'e de Paris-Sud in 1994-1995.  In this period, the original conjecture was fully proved, as we shall discuss later.   

In the next sections, after presenting the main definitions regarding regular Cantor sets and concepts of fractal geometry, we will present the above mentioned result by Jacob and Takens and then discuss the collaborations of Gugu and Jacob with Jean-Christophe related to homoclinic bifurcations and fractal geometry.

{\bf Acknowledgements:} We would like to warmly thank Carlos Matheus for his excellent suggestions for the present work, specially concerning the description of the long paper \cite{PY8} by Palis and Yoccoz. We are also grateful to the referee of this paper for his very valuable corrections and remarks. Finally, we would like to thank Amaury Alvarez and S\'ergio Vaz for helping in the editoring of the present text. 

\section{Homoclinic bifurcations in small dimensions: the work by Palis and Takens}

The first natural problem related to homoclinic bifurcations is the study of {\it homoclinic explosions} on surfaces: We consider one-parameter families $(\varphi_{\mu})$, $\mu \in (-1,1)$ of diffeomorphisms of a surface for which $\varphi_{\mu}$ is uniformly hyperbolic for $\mu<0$, and $\varphi_0$ presents a quadratic homoclinic tangency associated to a hyperbolic periodic point (which may belong to a {\it horseshoe} - a compact, locally maximal, hyperbolic invariant set of saddle type). It unfolds for $\mu>0$ creating locally two transverse intersections between the stable and unstable manifolds of (the continuation of) the periodic point. A main question is: what happens for (most) positive values of $\mu$? The following figure depicts such a situation for $\mu=0$.

%\begin{figure}
%\centering
\begin{center}
\includegraphics[width=5cm]{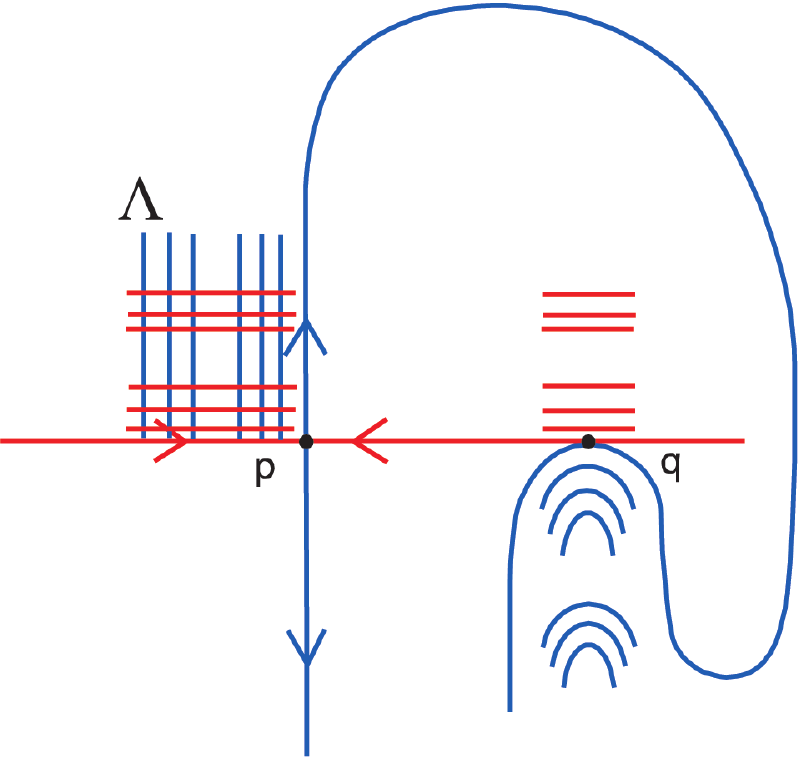}
\end{center}
%\end{figure}

Fractal sets appear naturally in Dynamical Systems and fractal dimensions when we try to measure fractals. They are essential to describe most of the main results in this presentation.

Given a metric space $X$, it is often true that the minimum number $N(r)$ of balls of radius $r$ needed to cover $X$ is roughly proportional to $1/r^d$, for some positive constant $d$, when $r$ becomes small. In this case, $d$ will be the {\it box dimension} of $X$. More precisely, we define the (upper) box dimension of $X$ as 
$$d(X)=\limsup_{r\to 0} \frac{\log N(r)}{-\log r}.$$

The notion of {\it Hausdorff dimension} of a set is more subtle, but more useful. The main difference with the notion of box dimension is that, while the box dimension is related to coverings of $X$ by small balls of equal radius, the Hausdorff dimension deals with arbitrary coverings of $X$ by balls of small (but not necessarily equal) radius. 

Given a countable covering $\cal U$ of $X$ by
balls, ${\cal U} = (B(x_i,r_i))_{i\in \N}$, we define its norm $||{\cal U}||$ as 
$||{\cal U}|| =
\max\{r_i,\,i\in \N\}$ (where $r_i$ is the radius of the ball $B(x_i,r_i)$). Given
$s\in\re_+$, we define ${\cal H}_s({\cal U}) =\sum_{i\in \N} r_i^s$.

The Hausdorff $s$-measure of $X$ is 
$${\cal H}_s(X) = \lim_{\epsilon\to
0} \inf_{{\cal U} \text{ covers } X \atop ||{\cal U}||<\epsilon}{\cal H}_s({\cal U}).$$

One can show that there is an unique real number, the {\it Hausdorff
dimension} of $X$, which we denote by $HD(X)$, such that $s < HD(X) \Rightarrow
{\cal H}_s(X) = +\infty$ and $s > HD(X) \Rightarrow {\cal H}_s(X) = 0$.

For ``well-behaved" sets $X$ - in particular for the sets which we study in this paper (regular Cantor sets and horseshoes), the box and Hausdorff dimensions of $X$ coincide.
    
Regular Cantor sets on the line play a fundamental role in dynamical systems and notably also in some problems in number theory. They are defined by expanding maps and have some kind of self-similarity property: small parts of them are diffeomorphic to big parts with uniformly bounded distortion (see a precise definition in the next section). In both settings, dynamics and number theory, a key question is whether the arithmetic difference (see definition below) of two such sets has non-empty interior.

A horseshoe $\Lambda$ in a surface is locally diffeomorphic to the Cartesian product of two regular Cantor sets: the so-called {\it stable} and {\it unstable} Cantor sets $K^s$ and $K^u$ of $\Lambda$, given by intersections of $\Lambda$ with local stable and unstable manifolds of some points of the horseshoe. The Hausdorff dimension of $\Lambda$, which is equal to the sum of the Hausdorff dimensions of $K^s$ and $K^u$, plays a fundamental role in several results on homoclinic bifurcations associated to $\Lambda$.

From the dynamics side, in the eighties, Jacob and Takens (\cite{PT}, \cite{PT1})
proved the following theorem about homoclinic bifurcations associated to a hyperbolic set:

\begin{theorem}
Let $(\varphi_{\mu})$, $\mu \in (-1,1)$ be a family of diffeomorphisms of a surface presenting a homoclinic explosion at $\mu=0$ associated to a periodic point belonging to a horseshoe $\Lambda$. Assume that $HD(\Lambda)<1$. Then
$$\lim_{\delta \to 0}\frac{m(H\cap [0,\delta])}{\delta}=1,$$
where $H:=\{\mu>0 \mid \varphi_{\mu}\,\hbox{  is uniformly hyperbolic}\}$.
\end{theorem}

A central fact used in the proof of this result is that if $K_1$ and $K_2$ are regular Cantor sets on the real line such that the sum of their Hausdorff dimensions is smaller than one, then $K_1-K_2 =
\{x-y \mid x \in K_1, y \in K_2\}=\{t\in \R|K_1\cap(K_2+t)\ne \emptyset\}$ (the {\it arithmetic difference} between $K_1$ and $K_2$) is a set of zero Lebesgue measure (indeed of Hausdorff dimension smaller than 1). On the occasion, looking for some kind of characterization property for this phenomenon, Jacob conjectured (see \cite{P}, \cite{P1}) that for generic pairs of regular Cantor sets $(K_1, K_2)$ of the real line either $K_1- K_2$ has zero measure or else it contains an interval (the last case should correspond in homoclinic bifurcations to open sets of tangencies). A slightly stronger statement is that,
 if $K_1$ and $K_2$ are generic regular Cantor sets and the sum of their Hausdorff dimensions is bigger than 1, then $K_1-K_2$ contains intervals. 

Another motivation for the conjecture was Newhouse's work in the seventies, when he introduced the concept of thickness of a regular Cantor set, another fractal invariant associated to Cantor sets on the real line. It was used in \cite{N1} to exhibit open sets of
diffeomorphisms with persistent homoclinic tangencies, therefore with no
hyperbolicity. It is possible (\cite{N2}) to prove that, under a dissipation hypothesis, in such an open set there is
a residual set of diffeomorphisms which present infinitely many coexisting
sinks. In \cite{N3}, it is proved that under generic hypotheses every family of
surface diffeomorphisms that unfold a homoclinic tangency goes through such an
open set.  It is to be noted that in the case described above with $HD(\Lambda)<1$ (as studied in \cite{PT1}) these sets have zero density.  See \cite{PT2} for a detailed presentation of these results. An important related question by Palis is whether the sets of parameter values corresponding to infinitely many coexisting sinks have typically zero Lebesgue measure.

An earlier and totally independent development had taken place in number theory. In 1947, M. Hall (\cite{H}) proved that any real number can be written as the sum of two numbers whose continued fractions  coefficients (of positive index) are bounded by $4$. More precisely, if $C(4)$ is the regular subset formed of such numbers in $[0,1]$, then one has $C(4)+C(4) = [\sqrt2 -1, 4(\sqrt 2 -1)]$.

In the next section we will discuss the positive solution of Palis' conjecture in the $C^k$-topology, $k>1$, by Gugu and Jean-Christophe, and some dynamical consequences of it on the study of homoclinic bifurcations. Before that, we will discuss two joint papers by Jacob and Jean-Christophe on this subject.

In the paper \cite{PY4}, they proved the following result: given a horseshoe $\Lambda$ with $HD(\Lambda)>1$ of a diffeomorphism $\varphi$, assume that $\varphi$ presents a quadratic tangency associated to a periodic point of $\Lambda$. Consider typical two-parameter families $\varphi_{s,t}$ of diffeomorphisms such that $\varphi_{0,0}=\varphi$ and, for each $s$ small, $\varphi_{s,0}$ presents a quadratic tangency associated to a periodic point (an analogous, in two parameters, to a homoclinic explosion). Then, for almost every $s$, the set of $t>0$ such that $\varphi_{s,t}$ presents a tangency between stable and unstable leaves of the continuation of $\Lambda$ has positive (upper) Lebesgue density at $t=0$, and so, contrary to the case of Palis-Takens' theorem, there is not full density of hyperbolicity at the initial bifurcation point. 

In the other paper \cite{PY5}, they proved that, densely, if $K_1$ and $K_2$ are regular Cantor sets such that $HD(K_1)+HD(K_2)>1$, then $K_1-K_2$ has non-empty interior. More precisely, if $HD(K_1)+HD(K_2)>1$, then there are $C^{\infty}$ diffeomorphisms $h$ of the real line arbitrarily close to the identity such that $K_1-h(K_2)$ has non-empty interior.

An important tool in the proof of both results is Marstrand's theorem, which we shall discuss later.
\section{Regular Cantor sets and homoclinic bifurcations}

A Cantor set $K\subset {\mathbb R}$ is a {\it $C^k$-regular Cantor set}, $k\ge 1$, if:

\begin{itemize}

\item[i)] there are disjoint compact intervals $I_1,I_2,\dots,I_r$ such that $K\subset I_1 \cup \cdots\cup I_r$ and the boundary of each $I_j$ is contained in $K$;

\item[ii)] there is a $C^k$ expanding map $\psi$ defined in a neighbourhood of $I_1\cup I_2\cup\cdots\cup I_r$ such that, for each $j$, $\psi(I_j)$ is the convex hull of a finite union of some of these intervals $I_s$. Moreover, we suppose that $\psi$ satisfies:

\begin{itemize}

\item[ii.1)] for each $j$, $1\le j\le r$ and $n$ sufficiently big, $\psi^n(K\cap I_j)=K$;

\item[ii.2)] $K=\bigcap\limits_{n\in\mathbb N} \psi^{-n}(I_1\cup I_2\cup\cdots\cup I_r)$.

\end{itemize}
\end{itemize}

\begin{remark}
If $k$ is not an integer, say $k=m+\alpha$, with $m\ge 1$ integer and $0<\alpha<1$ we assume that $\psi$ is $C^m$ and $\psi^{(m)}$ is $\alpha$-H\"older.
\end{remark}

We say that $\{I_1,I_2,\dots,I_r\}$ is a {\it Markov partition} for $K$ and that $K$ is {\it defined} by $\psi$.

\begin{remark}
In general, we say that a set $X \subset \R$ is a Cantor set if $X$ is compact, without isolated points and with empty interior. Cantor sets in $\R$ are homeomorphic to the classical ternary Cantor set $K_{1/3}$ of the elements of $[0,1]$ which can be written in base 3 using only digits $0$ and $2$. The set $K_{1/3}$ is itself a regular Cantor set, defined by the map $\psi:[0,1/3]\cup [2/3,1] \to \R$ given by $\psi[x]=3x-\lfloor 3x \rfloor$.
\end{remark}

An {\it interval of the construction of the regular Cantor set $K$} is a connected component of $\psi^{-n}(I_j)$ for some $n\in\mathbb N$, $j\le r$.

Given $s \in [1,k]$ and another regular Cantor set $\tilde K$, we say that $\tilde K$ is close to $K$ in the $C^s$ topology if $\tilde K$ has a Markov partition $\{\tilde I_1,\tilde I_2,\dots,\tilde I_r\}$ such that the interval $\tilde I_j$ has endpoints close to the endpoints of $I_j$, for $1 \le j \le r$ and $\tilde K$ is defined by a $C^s$ map $\tilde \psi$ which is close to $\psi$ in the $C^s$ topology.

The $C^{1+}$-topology is such that a sequence $\psi_n$ converges to $\psi$ if there is some $\alpha>0$ such that $\psi_n$ is $C^{1+\alpha}$ for every $n\ge 1$ and $\psi_n$ converges to $\psi$ in the $C^{1+\alpha}$-topology.

The concept of stable intersection of two regular
Cantor sets was introduced in \cite{M}: two Cantor sets $K_1$ and $K_2$ have stable intersection if there is a neighbourhood $V$ of $(K_1,K_2)$ in the set of pairs of $C^{1+}$-regular Cantor sets such that $(\widetilde
K_1, \widetilde K_2) \in V \Rightarrow \widetilde K_1 \cap \widetilde K_2 \ne \emptyset$. 

In the same paper conditions based on
renormalizations were introduced to ensure stable intersections, and applications of stable intersections to homoclinic bifurcations were obtained: roughly speaking, if some translations of the stable and unstable regular Cantor sets associated to the horseshoe at the initial bifurcation parameter $\mu=0$ have stable intersection then the set $\{\mu>0 \mid \varphi_{\mu}\,\,\hbox{ presents persistent homoclinic tangencies}\}$ has positive Lebesgue density at $\mu=0$. It was also shown that this last phenomenon can coexist with positive density of hyperbolicity in a
persistent way. 

Besides, the following question was posed in \cite{M}:
Does there exist a dense (and automatically open) subset ${\cal U}$ of
$$\Omega^{\infty} =\{(K_1,K_2); K_1, K_2\,\,\hbox{are} \,C^{\infty}-\hbox{regular Cantor sets and}\, HD(K_1) + HD(K_2) > 1\}$$
%$$\Omega^{\infty} =\{(K_1,K_2); K_1, K_2\,\,\,\, C^{\infty}-\hbox{regular Cantor sets}\mid HD(K_1) + HD(K_2) > 1\}$$
such that $(K_1,K_2) \in {\cal U} \Rightarrow \exists\, t\in\re$ such that $(K_1,K_2+t)$ has stable intersection?
A positive answer to this question implies a strong version of Palis'
conjecture. Indeed, $K_1-K_2=\{t\in \R\mid K_1\cap (K_2+t)\ne \emptyset\}$, so, if $(K_1,K_2+t)$ has stable intersection then $t$ belongs persistently to the interior of $K_1-K_2$.

The results of \cite{MY} gave an affirmative answer to this question, proving the following

\begin{theorem}
There is an open and dense set \,\,${\cal U} \subset \Omega^{\infty}$ such that if $(K_1,K_2) \in \cal{U}$, then $I_s(K_1,K_2)$ is dense in $K_1-K_2$ and $HD((K_1-K_2)\backslash I_s(K_1,K_2)) < 1$, where 
$$I_s(K_1,K_2) := \{t \in \mathbb{R} \mid K_1\hbox{ and }(K_2+t)\hbox{ have stable intersection}\}.$$
\end{theorem}

The same result works if we replace stable intersection
by $d$-stable intersection, which is defined by asking that any pair $(\widetilde K_1, \widetilde K_2)$ in some neighbourhood of $(K_1,K_2)$  satisfies $HD(\widetilde K_1 \cap \widetilde K_2) \ge d$: most pairs of Cantor sets $(K_1,K_2) \in {\Omega}^{\infty}$ have $d$-stable intersection for any $d < HD(K_1) + HD(K_2)-1$.

The open set $\cal U$ mentioned in the above theorem is very large in ${\Omega}^\infty$ in the sense that generic $n$-parameter families in ${\Omega}^\infty$ are actually contained in $\cal U$.

The proof of this theorem depends on a sufficient condition for the existence of stable intersections of two Cantor sets, related to a notion of renormalization, based on the fact that small parts of regular Cantor sets are diffeomorphic to the whole set: the existence of a {\it recurrent compact set} of relative positions of limit geometries of them. Roughly speaking, it is a compact set of relative positions of regular Cantor sets such that, for any relative position in such a set, there is a pair of (small) intervals of the construction of the Cantor sets such that the renormalizations of the Cantor sets associated to these intervals belong to the interior of the same compact set of relative positions. 

The main result is reduced to prove the existence of recurrent compact sets of relative positions for most pairs of regular Cantor sets whose sum of Hausdorff dimensions is larger than one. A central argument in the proof of this fact is a probabilistic argument \`a la Erd\H os: we construct a family of perturbations with a large number of parameters and show the existence of such a compact recurrent set with large probability in the parameter space (without exhibiting a specific perturbation which works). 

Another important ingredient in the proof is the Scale Recurrence Lemma, which, under mild conditions on the Cantor sets (namely that at least one of them is not essentially affine), there is a recurrent compact set for renormalizations at the level of relative scales of limit geometries of the Cantor sets. This lemma is the fundamental tool in the paper \cite{M2}, in which it is proved that, under the same hypothesis above, if $K$ and $K'$ are $C^2$-regular Cantor sets, then $HD(K+K')=\min\{1, HD(K)+HD(K')\}$.

These results have, somewhat surprisingly, deep consequences in number theory, which we discuss below.

Let $\al$ be an irrational number. According to Dirichlet's theorem, the
inequality \linebreak $|\al-\frac pq|<\frac1{q^2}$ has infinitely many
rational solutions $\frac pq$. Hurwitz  and Markov improved this result by proving that
$|\al-\frac pq|<\frac1{\sqrt 5 q^2}$ also has infinitely many rational
solutions $\frac pq$ for any irrational $\al$, and that $\sqrt 5$ is the largest
constant that works for any irrational $\al$. However, for particular values of
$\al$ we can improve this constant.

More precisely, we define $k(\al):=\sup\{k>0\mid  |\al-\frac pq|<\frac 1{
kq^2}$ has infinitely many rational solutions $\frac pq\}=
\limsup_{p,q\to\infty}\, |q(q \al-p)|^{-1}$. We have
$k(\al)\ge\sqrt 5$, $\forall \al\in\re\setminus {\mathbb Q}$ and $k\left(
\frac{1+\sqrt 5}2 \right)=\sqrt 5$. We will consider the set $L=\{k(\al) \mid
\al\in\re\setminus {\mathbb Q}$, $k(\al)<+\infty\}$.

This set is called the {\it Lagrange spectrum}. Hurwitz-Markov theorem determines the smallest
element of $L$, which is $\sqrt 5$. This set $L$ encodes many diophantine properties of real
numbers. The study of the geometric structure of $L$ is a classical subject.

The classical {\it Markov spectrum} also has an arithmetic definition, related to minima in the integers of indefinite binary quadratic forms.

The above mentioned result by M. Hall implies that the Markov and Lagrange spectra contain the whole half-line $[6,+\infty)$. The main result of \cite{M3} states that the Hausdorff dimensions of the intersections of the Markov and Lagrange spectra with half-lines of the type $(-\infty, t)$ always coincide, and assume all values in the interval $[0,1]$ when $t$ varies in the interval $[3,\sqrt{12})$. 

There are many generalizations of the Markov and Lagrange spectra in geometry and dynamics (in particular, the classical Markov and Lagrange spectra are respectively the sets of maximum heights and asymptotic maximum heights of geodesics in the modular surface - a fact explained by Jean-Christophe to Gugu). The ideas and results of \cite{MY} are important tools also in several recent works on geometric properties of these generalized spectra (see \cite{CMM}, \cite{Ma}, \cite{MR1}, \cite{MR2} and \cite{MPR}).

An important result in fractal geometry which is used in \cite{MY} is the famous Marstrand's theorem (\cite{Mar}), according to which, given a Borel set $X \subset {\mathbb R}^2$ with $HD(X)>1$ then, for almost every $\lambda \in {\mathbb R}$, $\pi_{\lambda}(X)$ has positive Lebesgue measure, where $\pi_{\lambda}:{\mathbb R}^2 \to {\mathbb R}$ is given by $\pi_{\lambda}(x,y)=x-\lambda y$. In particular, if $K_1$ and $K_2$ are regular Cantor sets with $HD(K_1)+HD(K_2)>1$ then, for almost every $\lambda \in {\mathbb R}$, $K_1-\lambda K_2$ has positive Lebesgue measure. Gugu and Yuri Lima gave combinatorial alternative proofs of Marstrand's theorem, first in the case of Cartesian products regular Cantor sets (\cite{LiM1}) and then in the general case (\cite{LiM2}).

In \cite{MY2}, Gugu and Jean-Christophe proved the following fact concerning
generic homoclinic bifurcations associated to two dimensional 
saddle-type hyperbolic sets (horseshoes) with Hausdorff dimension bigger than one: typically there are translations of the stable and unstable Cantor sets having stable intersection, and so it yields open sets of stable tangencies in the parameter line with positive density at the initial bifurcation value. Moreover, the union of such a set with the hyperbolicity set in the parameter line generically has full density at the initial bifurcation value. This extends the results of \cite{PY4}.

\section{Non-uniformly hyperbolic horseshoes}
In the most recent work by Jacob and Jean-Christophe (\cite{PY8}; see also \cite{PY6} and \cite{PY7}), they propose to advance considerably the current knowledge on the topic of bifurcations of heteroclinic cycles for smooth ($C^{\infty}$) parametrized families $\{g_t, t\in \R\}$ of surface diffeomorphisms. They assume that $g_t$ is hyperbolic for $t<0$ and $|t|$ small and that a quadratic tangency $q$ is formed at $t=0$ between the stable and unstable lines of two periodic points, not belonging to the same orbit, of a (uniformly hyperbolic) horseshoe $\Lambda$ and that such lines cross each other with positive relative speed as the parameter evolves, starting at $t=0$ near the point $q$. They also assume that, in some neighbourhood $W$ of the union of $\Lambda$ with the orbit of tangency $O(q)$, the maximal invariant set for $g_0=g_{t=0}$ is $\Lambda \cup O(q)$, where $O(q)$ denotes the orbit of $q$ for $g_0$.

%\begin{figure}
%\centering
\begin{center}
\includegraphics[width=12cm]{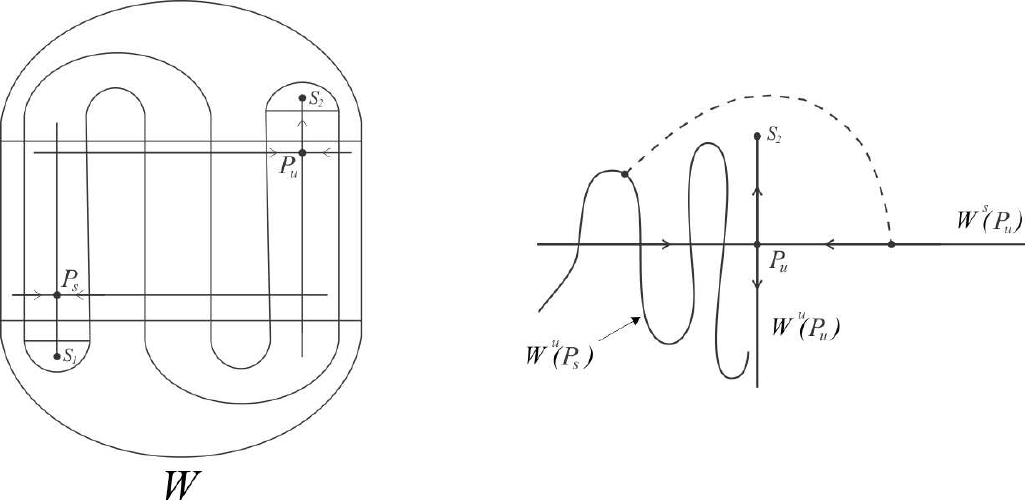}
\end{center}
%\end{figure}

A main novelty is that they allowed the Hausdorff dimension $HD(\Lambda)$ to be bigger than one: they assume $HD(\Lambda)$ to be bigger than one, but not much bigger (more precisely, they assume that, if $d_s$ and $d_u$ are the Hausdorff dimensions of, respectively, the stable and unstable Cantor sets of $g_0$ then $(d_s+d_u)^2+\max(d_s, d_u)^2<d_s+d_u+\max(d_s, d_u)$). Then, for most small values of $t$, $g_t$ is a ``non-uniformly hyperbolic" horseshoe in $W$, and so $g_t$ has no attractors nor repellors in $W$. Most small values of $t$, and thus most $g_t$, here means that $t$ is taken in a set of parameter values with Lebesgue density one at $t=0$.   

The construction of non-uniformly hyperbolic horseshoes for most parameters is a highly non-trivial counterpart of Jean-Christophe's proof \cite{Y-J} (based on the so-called \emph{Yoccoz puzzles}) of the celebrated Jakobson's theorem in the context of heteroclinic explosions. 

More concretely, the maximal invariant set $W$ can be described in terms of an iterated\footnote{It is worth to note that we do not specify whether we are iterating in the past or future: in fact, the past and future iterations are treated on equal foot because the arguments of Jacob and Jean-Christophe are \emph{time-symmetric} (and this is of course in sharp contrast with the arguments of Benedicks-Carleson \cite{BC} for the construction of H\'enon attractors).} system of functions $\{g_{(a),t}\}_{a\in\mathcal{A}}$ related to the (uniformly hyperbolic) horseshoe $\Lambda$ and a \emph{folding map} $G_t$ related to the orbit $O(q)$ of heteroclinic tangency. In this context, the idea is to define the \emph{candidate} for non-uniformly hyperbolic horseshoe $K\subset W$ as the set of points captured by nested sequences of domains of \emph{affine-like} hyperbolic maps (in the sense of \cite{PY6}) given by \emph{certain} compositions of $g_{(a),t}$ and $G_t$ with \emph{hope} that the set $\mathcal{E}=W\setminus K$ of points not captured by this scheme is <<exceptional>> (or <<small>>) for most choices of small values of $t$, say $t\in [\varepsilon_0, 2\varepsilon_0]$ with $\varepsilon_0\ll 1$.  

For an \emph{arbitrary} parameter $t\in [\varepsilon_0, 2\varepsilon_0]$, there is \emph{no} reason for the candidate $K$ to display non-uniform hyperbolicity. Nevertheless, Jacob and Jean-Christophe show that $K$ is non-uniformly hyperbolic for most $t\in[\varepsilon_0,2\varepsilon_0]$ via the following \emph{colossal}\footnote{The original article \cite{PY8} takes 217 pages to justify this whole process.} inductive parameter selection process. 

For each parameter interval $I\subset [\varepsilon_0, 2\varepsilon_0]$, Jacob and Jean-Christophe construct a \emph{canonical} class $\mathcal{R}(I)$ of affine-like maps of certain ($I$-persistent) simple and parabolic compositions of $g_{(a),t}$ and $G_t$. The non-uniform hyperbolicity of the portion of $K$ captured by the domains of the affine-like maps in the class $\mathcal{R}(I)$ can be proved whenever $I$ is \emph{strongly regular}: very roughly speaking, this means that $\mathcal{R}(I)$ contains few \emph{bicritical} elements corresponding to returns of the almost tangency region to itself (or, equivalently, to <<almost tangencies of higher orders>>). In particular, it is not surprising that the candidate $K$ is a non-uniformly hyperbolic horseshoe\footnote{Furthermore, Jacob and Jean-Christophe show that the set $\mathcal{E}=W\setminus K$ is <<exceptional>> in many senses.} whenever $t$ is a \emph{strongly regular parameter}, that is, a parameter $\{t\}=\bigcap\limits_{n\in\mathbb{N}} I_n$ given by a nested sequence $I_{n+1}\subset I_n$ of strongly regular intervals. 

At this point, Jacob and Jean-Christophe want to complete the proof of their result by showing that most parameters $t\in [\varepsilon_0, 2\varepsilon_0]$ are strongly regular. For this sake, they fix an adequate parameter $0<\varepsilon_0\ll \tau\ll 1$, they define $\varepsilon_{n+1}:=\varepsilon_n^{1+\tau}$, and they employ the following induction scheme: 
\begin{itemize}
\item one starts by proving\footnote{This is the only place where the fact that $O(q)$ is associated to a \emph{heteroclinic} tangency is used.} that the initial interval $I_0=[\varepsilon_0,2\varepsilon_0]$ is strongly regular; 
\item at the $k$th step, they divide all strongly regular intervals $I_{k-1}$ of the previous step into 
$\lfloor\varepsilon_k^{-\tau}\rfloor$ \emph{candidate} intervals $I_k$ of sizes $\varepsilon_k$; 
\item they keep the strongly regular intervals $I_k$ and they discard the candidates failing the strong regularity test; 
\item they prove that a proportion $\leq \varepsilon_0^{\tau^2}$ of the candidates are discarded; 
\end{itemize}
Here, let us mention that the proof of the last item above relies on the interplay\footnote{This is an incarnation of a long tradition in Dynamical Systems: as Adrien Douady used to say in the context of complex quadratic polynomials and the Mandelbrot set, one has to <<plough in parameter space, and harvest in phase space>> \cite{Dou}.} between phase space and parameter space: more precisely, Jacob and Jean-Christophe show that the strong regularity requirements in phase space have a counterpart in parameter space thanks to the way certain dynamical objects move with $t\in[\varepsilon_0,2\varepsilon_0]$, and this is the main reason for the strong regularity of most values of $t$. 

Anyhow, the conclusion of this induction scheme is that the Lebesgue measure of the subset of strongly regular parameters $t\in [\varepsilon_0, 2\varepsilon_0]$ is $\geq \varepsilon_0(1-\varepsilon_0^{\tau^2})$. Since $0<\varepsilon_0\ll \tau\ll 1$, this completes our sketch of proof of the main result in the paper \cite{PY8}. 

Of course, Jacob and Jean-Christophe do not consider their result as the end of the line. Indeed, they expected the same results to be true for all cases $0<HD(\Lambda)<2$. However, to achieve that, it seems that their methods need to be considerably sharpened: it would be necessary to study deeper the dynamical recurrence of points near tangencies of higher order (cubic, quartic...) between stable and unstable curves. They also expected their results to be true in higher dimensions.

Finally, they hoped that the ideas introduced in that work might be useful in broader contexts. In the horizon lies a famous question concerning the standard family of area preserving maps (which is the family $(f_{\lambda})_{\lambda \in \R}$ of diffeomorphisms of the torus ${\mathbb T}^2=\R^2/{\mathbb Z}^2$ given by $f_{\lambda}(x, y)=(-y+2x+\lambda\sin(2\pi x), x) \pmod 1$): can we find sets of positive  Lebesgue measure in the parameter space such that the corresponding maps display non-zero Lyapunov exponents in sets of positive Lebesgue probability in phase space?

\newpage

\section{Words by Jean-Christophe and photos}

\hskip .5in A generous text by Jean-Christophe in the visitors book of IMPA's library, 1989:

%\begin{figure}
%\centering
\begin{center}
\includegraphics[width=14.5cm]{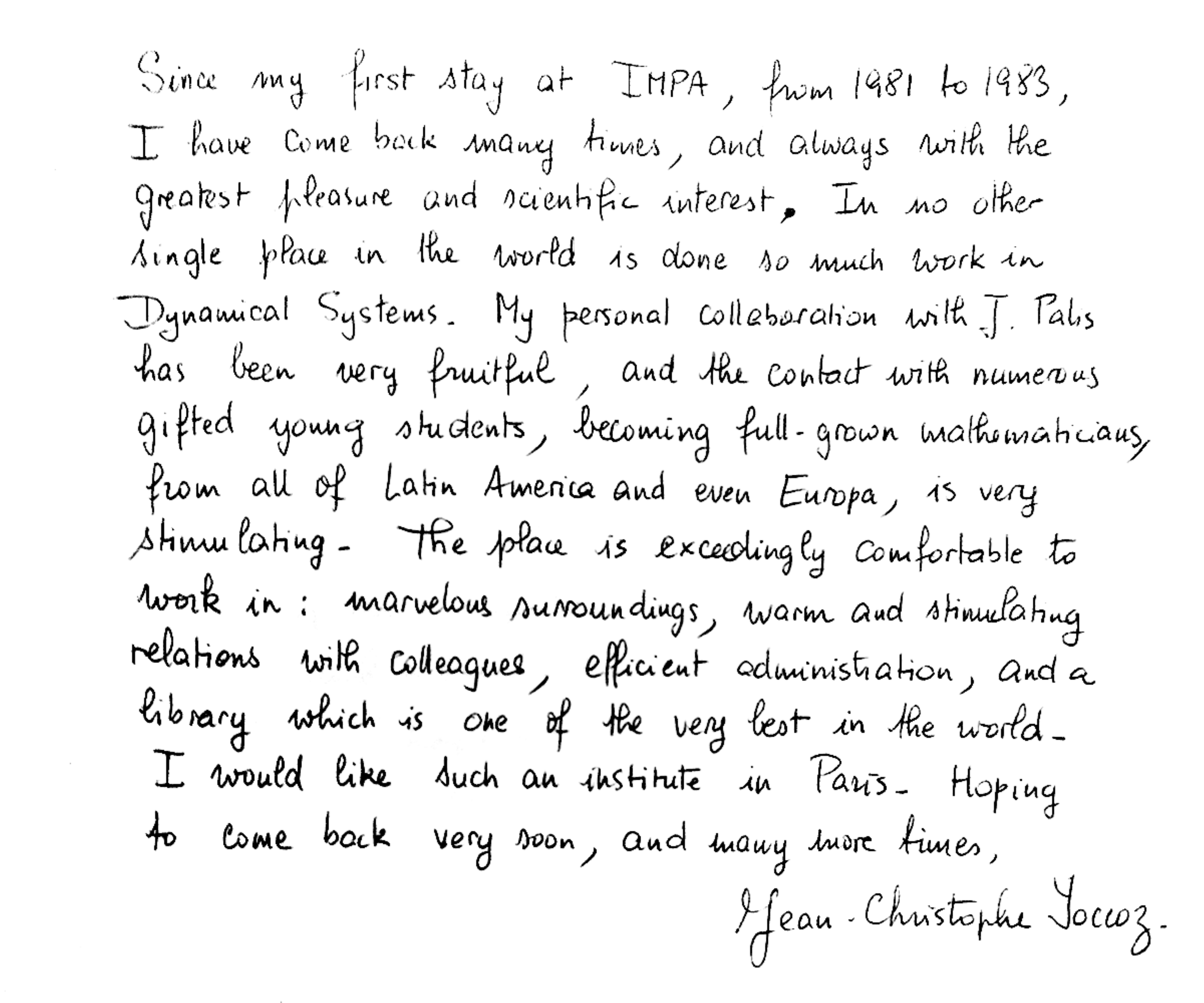}
\end{center}
%\end{figure}

Mauricio Peixoto, Gugu, Jean-Christophe and Elon Lima at IMPA, shortly after the Post-Doc of Gugu with Jean-Christophe:

%\begin{figure}
%\centering
\begin{center}
\includegraphics[width=10cm]{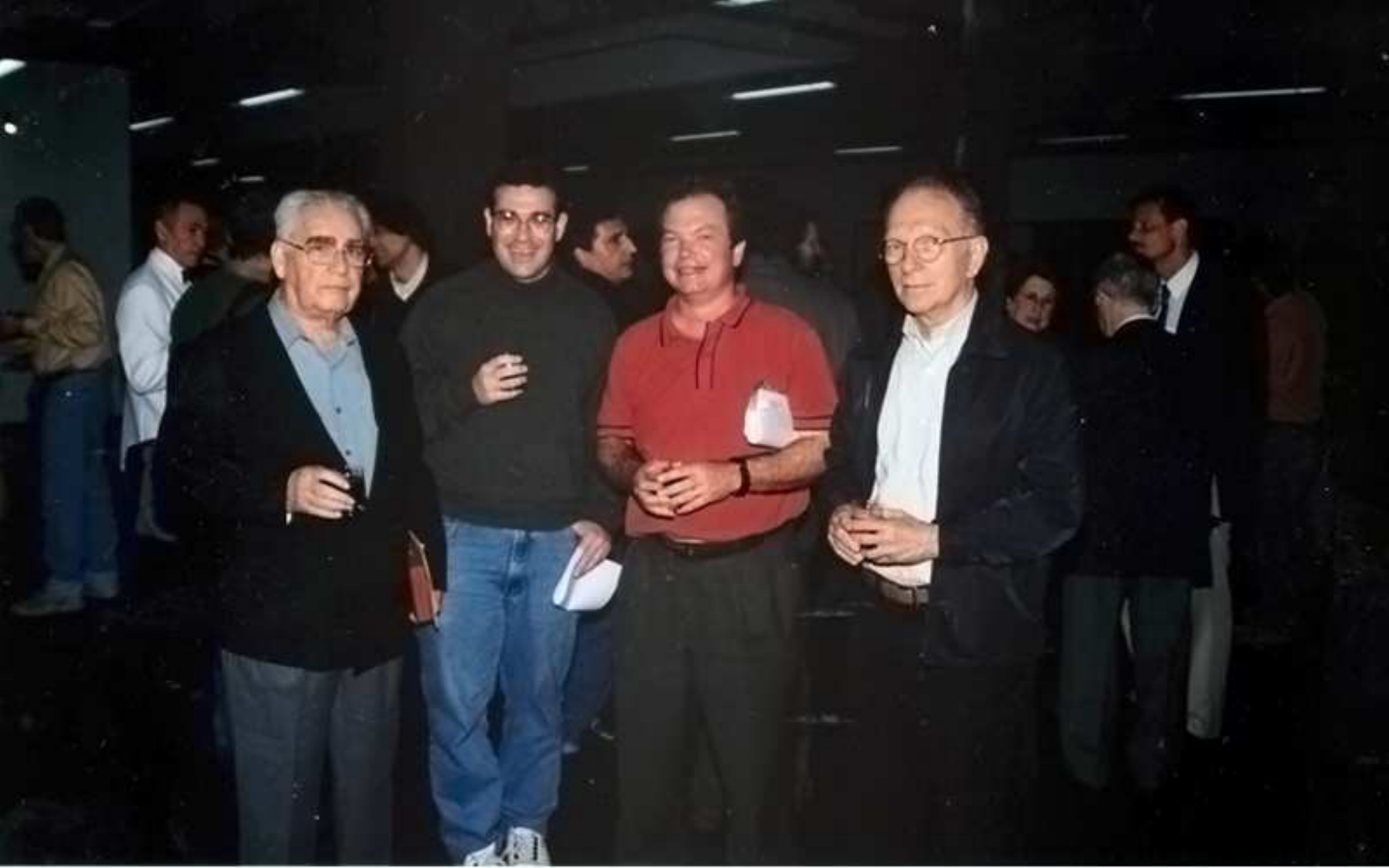}
\end{center}
%\end{figure}
%\newpage
\hskip 1in  A picture of Jacob and Jean-Christophe in San Marino:

%\begin{figure}
%\centering
\begin{center}
\includegraphics[width=9cm]{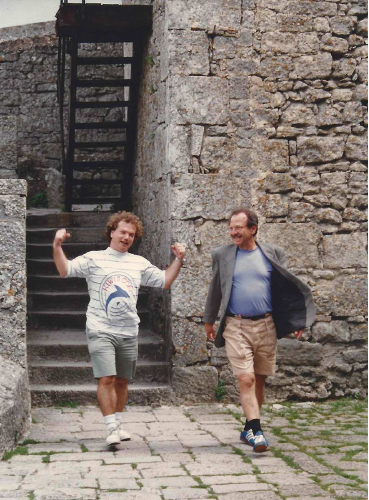}
\end{center}
%\end{figure}
\newpage 
{\bf Jacob Palis} 

IMPA - Estrada D. Castorina, 110 - 22460-320 - Rio de Janeiro - RJ - Brasil 

jpalis@impa.br

Jacob Palis is an Emeritus Researcher at IMPA and a member of its Dynamical Systems group since 1968. He has advised 42 Ph.D. students and has more than 240 Ph.D. descendants. His research covers several topics of Dynamical Systems, including global stability, bifurcations and fractal geometry. He posed many conjectures in the field, including a global one on the finitude of attractors for typical dynamical systems.
\vskip .15in
{\bf Carlos Gustavo Moreira}

IMPA - Estrada D. Castorina, 110 - 22460-320 - Rio de Janeiro - RJ - Brasil 

gugu@impa.br

Carlos Gustavo Moreira is a Researcher at IMPA and a member of its Dynamical Systems group since 1995. His research covers several branches of Dynamical Systems, including one-dimensional dynamics, bifurcations and fractal geometry. He also works on Combinatorics (including extremal and probabilistic combinatorics) and Number Theory (specially diophantine approximations).

%\frenchspacing

\end{document}